\documentclass[conference]{IEEEtran}
\IEEEoverridecommandlockouts
\usepackage[colorinlistoftodos]{todonotes}
\usepackage{booktabs}
\usepackage{cite}
\usepackage{amsmath,amssymb,amsfonts}
\usepackage{algpseudocode}
\usepackage{graphicx}
\usepackage{textcomp}
\usepackage{xcolor}

\def\BibTeX{{\rm B\kern-.05em{\sc i\kern-.025em b}\kern-.08em
    T\kern-.1667em\lower.7ex\hbox{E}\kern-.125emX}}
\usepackage[linesnumbered,ruled]{algorithm2e}
\SetKwIF{If}{ElseIf}{Else}{if}{}{else if}{else}{end if}%
\SetKwFor{While}{while}{}{end while}%
\SetKwRepeat{Do}{do}{while}
\SetKw{KwGoTo}{go to}
\begin{document}

\title{Heuristic for Optimisation of Dark Store Facility Locations for Quick Commerce Businesses \\

}
\author{\IEEEauthorblockN{Prithvi Dinesh Kewalramani}
\IEEEauthorblockA{
\textit{Indian Institute of Technology Bombay}\\
 190010056@iitb.ac.in\\
 prithvik1969@gmail.com}
\and
\IEEEauthorblockN{Harshad Khadilkar}
\IEEEauthorblockA{
\textit{Indian Institute of Technology Bombay}\\
 harshadk@iitb.ac.in\\
harshad.khadilkar@tcs.com }

}

\maketitle

\begin{abstract}
We present a fast, flexible heuristic for setting up warehouse locations for quick commerce businesses, with the goal of serving the largest number of customers under the constraints of delivery radius and maximum daily deliveries per warehouse. Quick commerce or direct-to-customer delivery businesses guarantee delivery within a specified time. Using experiments on various scenarios, we show that the proposed algorithm is flexible enough to handle variations such as non-uniform population distributions, variable travel times, and selection of multiple warehouse locations.
\end{abstract}


\section{Introduction}

Many quick-delivery apps guarantee delivery directly to customers within a short time, usually less than 30 minutes. These apps employ warehouses, called dark stores \cite{BULDEORAI2023101333}, to serve different areas in the city, with urban or suburban centre having multiple dark stores that serve customers within a service radius. It is an obvious business objective to open dark stores at locations where they can cater to the highest potential number of customers, thus maximising profits. This project aims to provide an efficient heuristic solution to the problem of finding the best locations for warehouses in a given population region. These businesses do not seek to service every potential customer within the region. Instead, the aim is to fulfil the delivery time promise to those customers who are eligible (i.e. are within a short radius of the fulfilment centres) while maximising the number of customers that fall in this category. This problem is known to be NP-hard \cite{Papadimitriou1978}, and we use a variant of the vehicle routing problem called the Multiple Travelling Salesman Problem (mTSP) with a single depot with time constraints. There have yet to be any research experiments on this problem, even though \cite{7963559}, \cite{Xu2018}, \cite{ahmed} discuss solutions for the mTSP problem with a single depot, which we will discuss in further sections.

Our system includes a set of customer locations on a map where we wish to plot an optimal location for our warehouses, one by one. Since the problem is NP-Hard, we design a heuristic for optimising the warehouse locations, with the well-known 2-opt algorithm \cite{goetschalckx1989vehicle} as an inner-loop optimiser. This substantially reduces the time-complexity of the approach, while our experiments show that the resulting solutions are still efficient. Modifications such as traffic settings and population distribution have also been tested. In the rest of this paper, we include a review of related literature in Section \ref{sec:related}, insights into the problem model and the approach in Section \ref{sec:problem}, and results and experiments in Section \ref{sec:results}.

We believe that the contribution of this paper is to present a simple, easy-to-use, efficient, scalable heuristic for designing a delivery network without the high sophistication required for traditional optimisation approaches. We have introduced an algorithm uniquely suited to solving the challenge of optimizing for warehouse locations based on delivery parameters including service radii, maximum trip time and daily delivery constraints. Our work also includes the integration of precise geographic and traffic considerations into our algorithm, offering a practical and relevant solution for urban logistics. These variables, given our time-constrained model, distinguish it from traditional problems in this domain.


\section{Related Literature} \label{sec:related}

The problem addressed in this project broadly derives from two sub-problems, namely the facility location problem (FLP) and the multiple travelling salesman problem (mTSP). Past literature on both topics has been reviewed in this section, in the context of how it is related to our work.

\subsection{Facility location problem}

The following studies address the problem of warehouse location and the optimisation of delivery parameters, to achieve their respective objectives.

\cite{XIFENG201345} work on a multi-objective optimization model that considers economic, service and environmental objectives for logistics network design. The model is based on the uncapacitated facility location problem (UFLP) with binary location variables and multiple objective functions, working with a greedy heuristic. 

\cite{su15108098} addresses the challenge of optimising last-mile delivery in e-commerce, with a focus on sustainability. It employs a two-phase approach that first uses a multi-criteria decision-making method to determine optimal pickup and delivery points and then uses a mathematical model for efficient vehicle routing with various service options. Their approach reduces costs, minimises environmental impact, and enhances customer satisfaction by determining appropriate trips based on service options.

\cite{KLOSE20054} provide a survey of various approaches to facility location models for selecting optimal locations, to serve the end customer effectively or to get a more efficient supply-chain model. The paper talks about two models relevant to our study, which include network location models and mixed-integer programming models.

\begin{enumerate}
\item Network location models models represent potential facilities and customers as a graph, with distances measured as shortest paths. They seek to minimize the total cost associated with customer-facility assignments ($x_{ij}$), taking into account the distance ($d_{ij}$) from each customer to a facility, with constraints that each customer is assigned to one facility, and each facility serves one customer. 

\item Mixed-integer programming models use a discrete set of potential facilities and integer variables to represent location and allocation decisions. The formulation constraints such as capacity, demand, cost, service, and routing. It minimizes the overall cost, a combination of facility opening costs ($f_i$) and transportation costs ($c_{ij}$), considering the demand of each customer ($d_j$) and the capacity of each facility ($Q$). 
\end{enumerate}

Our work is on a variant of FLP, but in our case, we do not actually serve all of the customers on the map, but only aim to pick the location(s) that can serve the maximum potential customers feasibly within the given constraints of service radius, delivery time, and volume of daily deliveries. Hence it requires a different kind of optimisation model altogether, while also considering the factors of variable travel times and population densities.

\subsection{Multiple travelling salesman problem}

The multiple travelling salesman problem is a well known NP-hard problem, which has a multitude of applications in logistics and deliveries. Most attempts to provide solutions have been heuristics, some of which include the genetic algorithm \cite{6721865}, maximum entropy based algorithm \cite{7963559}, a two-phase heuristic \cite{Xu2018} and the ant colony optimization algorithm \cite{1691778}.

A proposed exact algorithm has been suggested by Ahmed et al. \cite{ahmed}, which involves a modified Lexisearch algorithm. It organises feasible solutions into blocks led by block leaders, then calculates lower bounds for these blocks, comparing them to the current best solution. When no superior solution is found in a block, the algorithm proceeds to the next. This iterative process continues until an optimal solution is obtained. 

\cite{7963559} propose a maximum entropy principle based approach for a few variants of the mTSP, including mTSP with a single depot, which is relevant to our work. Their algorithm reinterprets TSP as a resource allocation challenge, utilising a FLP framework. It employs a Deterministic Algorithm (DA) to iteratively to optimise city-facility associations, guided by an increasing annealing parameter. By maximizing marginal distributions, DA identifies optimal facility locations. 

\cite{Xu2018} introduce a two phase heuristic algorithm (TPHA) to solve the mTSP. In the first phase, an improved version of the K-means algorithm is used to group visited cities based on their locations and specific capacity constraints. In the second phase, a route planning algorithm is designed using a genetic algorithm (GA) combined with the roulette wheel method and elitist strategy for selection. The proposed method minimizes overall traveling distance, and outperforms GA-based route planning algorithms in solving the mTSP. 

 \cite{Mitrovi2002TheMT} talk about the multiple travelling salesman problem with time windows (mTSPTW), where a fleet of vehicles has to serve a set of locations within specified time windows, with the objective to minimize the number of vehicles needed. The approach involves the use of two types of precedence graphs, which capture the order relation among the locations imposed by the travel times and the time windows. The method uses the minimum chain decomposition of these graphs to find a lower and an upper bound for the minimum number of vehicles. 

Our approach uses the mTSP problem in the second part of the solution, as a base for serving customers within a fixed radius. However, instead of optimising to serve all the customers with $m$ salesmen, we only serve as many customers as feasible with the constraint on total travel time in one trip (one trip may have multiple deliveries, and each customer is to be served only once). Thus, ours works as a mTSP with time constraints, which is markedly different from an mTSP with time windows, considering the fact that the latter involves fixed time windows where each location has to be served. Our solution is more relevant from the perspective of quick delivery apps, where our goal is simply to maximise the number of deliveries in an area from a given number of depots, given the constraints on maximum delivery time and number of deliverymen.



\section{Problem description and Proposed Method} \label{sec:problem}
	
This study addresses the problem of determining the optimal locations for $N$ number of dark stores within a given geographical area, such that maximum total number of customers can be served. Here $N$ is an externally specified parameter. We assume that the coordinates of customers are provided \textit{a priori}. Considering $M$ number of customers on the map, the location of each customer $i$ is represented as $\vec{C_i}(x_i, y_i)$. For a potential warehouse location on the map, represented by $\vec{W_j}(x_j, y_j)$, the decision variable $X_{ij}$  represents whether the customer will be served by warehouse $j$ or not. Our goal is to get a list of the locations of $N$ warehouses on the map, such that the maximum number of customers are served, as so:
\begin{equation}
\max \left(\sum_{i=1}^{N} \sum_{j=1}^{M} X_{ij}\right)
\label{eq:obj}
\end{equation}

The maximisation is subject to the following constraints: 

\begin{enumerate}
    \item Each warehouse will serve customers only within its service radius ($t_{max}$). The service radius is defined in terms of total time needed to reach the edge of the circle with the warehouse as its centre.
$$||\vec{C_i}(x_i, y_i)-\vec{W_j}(x_j, y_j)|| \leq t_{max},$$
, where $||\vec{A}-\vec{B}||$ represents the time taken to travel between two locations $\vec{A}$ and $\vec{B}$. This can be approximated as a constant multiple of Euclidean distance (in the uniform traffic case), and can be derived from a real-time map interface (in the case with variable traffic zones). The traffic zone cases will be talked about in Section III-B.

\item The total time for each trip from a depot cannot exceed $t_{max}$. One trip is defined as the delivery-person leaving their depot, to serve their allotted customers one by one, until the last customer. The total trip time is, thus, the time elapsed until their last delivery. It is essential to note that the time required for the delivery-person to return to the depot after the last delivery is not taken into account for this metric.

$$t_{trip} \leq t_{max}$$

\item Each customer is served by exactly one warehouse:
\begin{equation}
\sum_{i=1}^{N} X_{ij} = 1, \quad \forall j \in \{1, 2, \ldots, M\}
\end{equation}
\item The number of delivery-people at each warehouse is given by $m$, and the maximum daily trips allowed to each delivery-person is given by $n$. The number of daily trips per warehouse, given by $D_j$ should not exceed the  product of $m$ and $n$.

 \begin{equation}
     D_j \leq mn
 \end{equation}
\end{enumerate}

As a result, we receive the final optimal locations for the warehouses:
$$\vec{W_k}(x_k, y_k) \forall k \in 1,2,3,... N $$

\subsection{Algorithm description}

Consider a grid of arbitrary dimension, which spans the potential customer locations as well as the potential fulfilment centre locations. We maximise the scalability of our proposed algorithm by implementing a multi-resolution approach, where the warehouse locations are first picked on a coarse grid and then iteratively refined using finer local grid search.

The initial step involves partitioning the map into a 10 $\times$ 10 grid, designating each grid point $\vec{G_j}(X_j, Y_j)$ as a potential dark store location in an iterative process. At each location, the set of points within the maximum serviceable radius is defined. Subsequently, the algorithm employs a modified 2-opt method \cite{goetschalckx1989vehicle}, as shown in Algorithm \ref{alg:one} and Algorithm \ref{alg:two}, to solve the traveling salesman problem for one delivery person at a time. Notably, this variant of the 2-opt method ensures that the delivery personnel return to their starting point after a specified time interval of $t_{max}$ minutes, thus mimicking a real-world scenario. The algorithm terminates when further improvements to the total time become infeasible, resulting in the final path, covering all serviceable customers. This provides a solution for scenarios with unlimited deliveries. The top $D$ deliveries, serving the maximum number of customers, are selected from this solution, determining the maximum number of customers that can be served when the respective grid point serves as the dark store location.

\begin{algorithm}[t]
\caption{Finding  total deliveries potential from each dark store location}
\label{alg:one}
\begin{algorithmic}
\Require $C_i(x_i, y_i), D, t_{max}$

\For {$\vec{G_i}(X_i, Y_i)$}  {
$\forall C_i\Rightarrow time(\vec{G_i}-\vec{C_j}) \leq t_{max}$

Enter the initial order of the customers as initial path

Put this path through Algorithm 2, calculate total path time
    
    \While {improvement to path is still potential}{
        Swap the next to visit $C_i$ with each other
        
        Set this as the new route
        
        Put this path through Algorithm 2, calculate total path time
    }
    Divide the final path into individual orders and pick the top $D$ orders with highest deliveries

    \Return sum of deliveries and all customers served
    }
\end{algorithmic}
\end{algorithm}

\begin{algorithm}[t]
\caption{Dividing the path obtained from the 2-opt function into individual deliveries}
\label{alg:two}
\begin{algorithmic}

\Require The route obtained using the 2-opt implementation
The origin is set as the starting point of the journey

Set $L=0$
\While{ There are still places left to deliver within the serviceable radius} {
$L=L+time(\vec{P}-\vec{C_i})$ 
\If{$L>t_{max}$}{
Add the origin location as next point in the path
Set $L=0$
}
\Else{ 
Add $C_i$ as the next point in the path
}

}
Add the origin as the last point in the path
\Return final path
\end{algorithmic}
\end{algorithm}

Following this, the algorithm is executed for all grid points on the map. The top 5\% of grid points, ranked by the number of customers served, are retained. Subsequently, the algorithm is applied again to smaller regions around these top 5\% grid points. Each such region is divided into a 4 $\times$ 4 grid, with dimensions reduced to 1/5th of the original size. This process of selecting and subdividing regions continues iteratively until the required level of precision is achieved.

Once the location of the first warehouse is determined, the customers served by this warehouse are removed from the map, and the entire process is repeated for the subsequent warehouses until the desired number of dark stores is established.

\begin{figure}[h]
    \centering
    \includegraphics[width=0.7\linewidth]{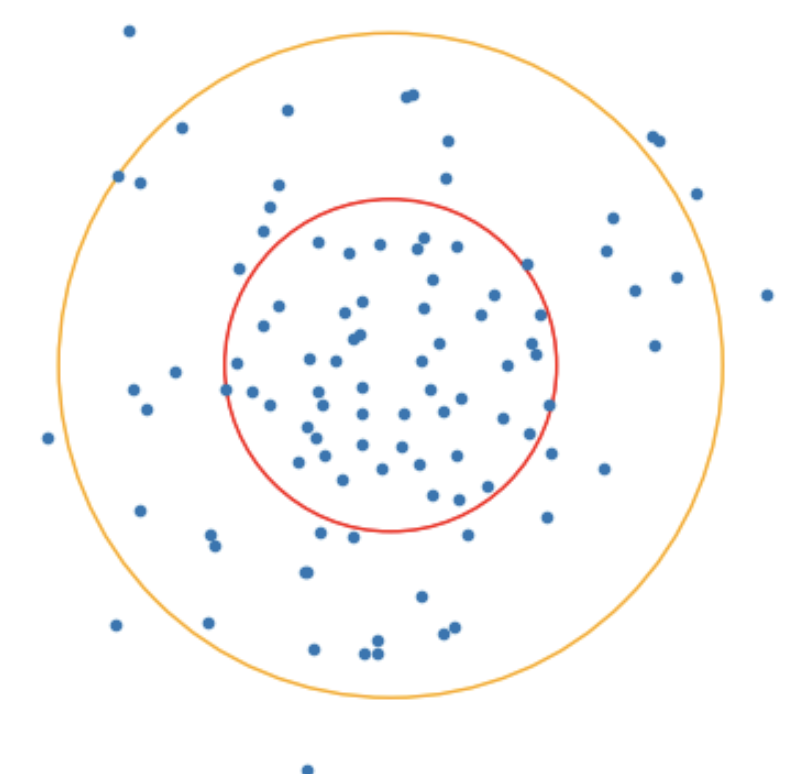}
   \caption{Traffic regions in the map, with the central region having the highest traffic multiplier, the middle one having a lower multiplier, and the outer region having the lowest multiplier. }
   \label{fig:traffic}
\end{figure}

Additionally, the algorithm can incorporate a uniform traffic scenario or variable traffic scenario. In the former, every unit distance in a map requires equal amount of time for traversal, while in the latter the map is classified into areas with varying traffic levels, such as high, moderate, and low, as shown in Fig. \ref{fig:traffic}. This implementation aligns with the approach employed in \cite{KOC201681}, where the time required to traverse a unit distance varies with how central the zone is in the given region.

\subsection{Comments on the approach}

Algorithms \ref{alg:one} and \ref{alg:two} aim to minimise the total path time for the delivery personnel, rather than directly optimising for the objective \eqref{eq:obj}. This is a deliberate design choice for facilitating the use of heuristics such as 2-opt, which accept perturbations to the solution only when there is an improvement to the objective. Incremental improvements are more difficult to achieve for discrete objectives like \eqref{eq:obj} than for continuous ones \cite{ouaarab2014discrete}. Therefore instead of directly pursuing the primary goal, the approach focuses on minimising the total path time to be covered, prior to selecting the top few routes. Essentially, the selection of an optimal final path capable of completing all deliveries within the shortest potential time serves as a proxy for the explicit maximization of total customers served within the first $N$ deliveries. 

\section{Results and Conclusions} \label{sec:results}

We have considered two population distribution models, namely Gaussian and uniform, and two types of traffic models, namely a uniform traffic model and a zone based traffic model. The distribution data comes from the Gaussian and random distribution functions, respectively, in the NumPy library in Python. In each of the cases, we compute the locations for 3 warehouses (which have been selected based on the maximum number of customers served in the iteration) with their routes for serving their respective customers. The accompanying results show these warehouse routes in the order of orange, green and red for the first, second, and third warehouses respectively. 

Below we consider the following 4 cases, and describe our observations.

\begin{figure}[h]
    \centering
    \includegraphics[width=0.7\linewidth]{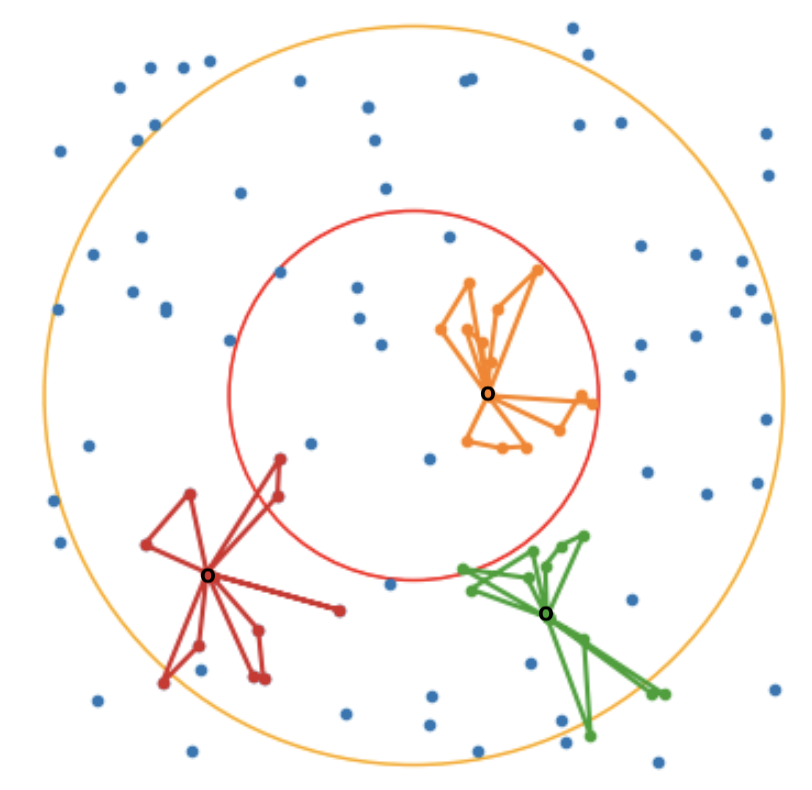}
   \caption{Warehouse distribution in the case where population distrubution and traffic are both uniform (Note that zones are still shown for comparison to other cases)}
   \label{fig:case4}
\end{figure}

\begin{figure}[h]
    \centering
    \includegraphics[width=0.7\linewidth]{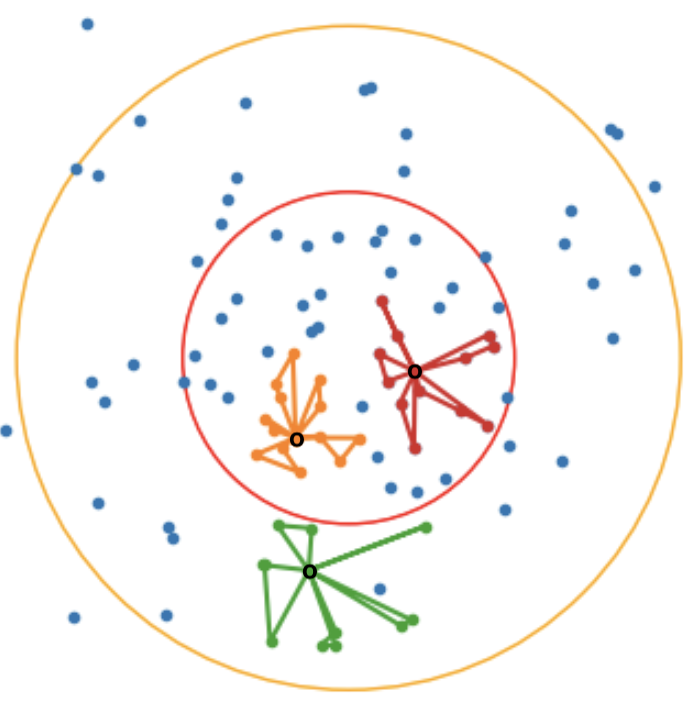}
   \caption{Warehouse distribution in the case where population distribution is gaussian, and traffic zones are present}
   \label{fig:case1}
\end{figure}

\begin{figure}[h]
    \centering
    \includegraphics[width=0.7\linewidth]{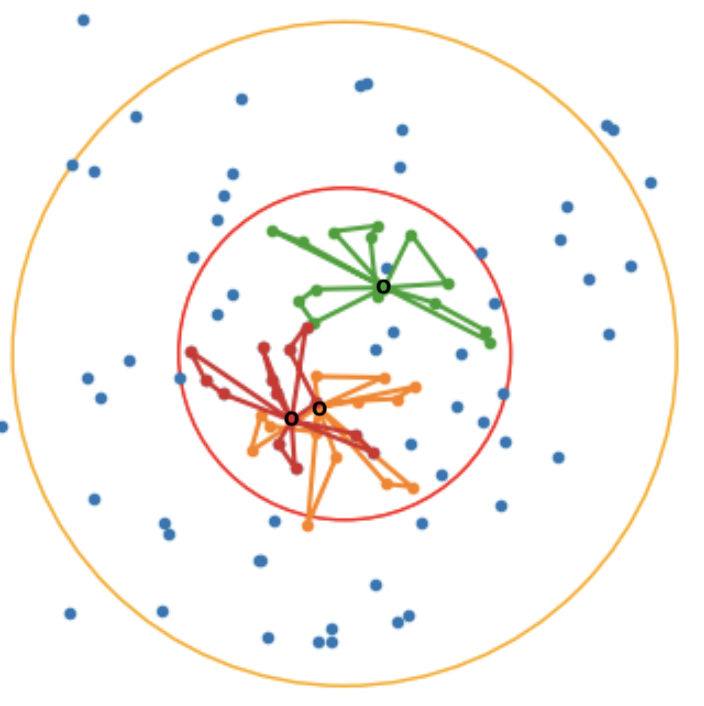}
   \caption{Warehouse distribution in the case where population distribution is gaussian, and traffic is uniform (Note that zones are still shown for comparison to other cases)}
      \label{fig:case2}
\end{figure}

\begin{figure}[h]
    \centering
    \includegraphics[width=0.7\linewidth]{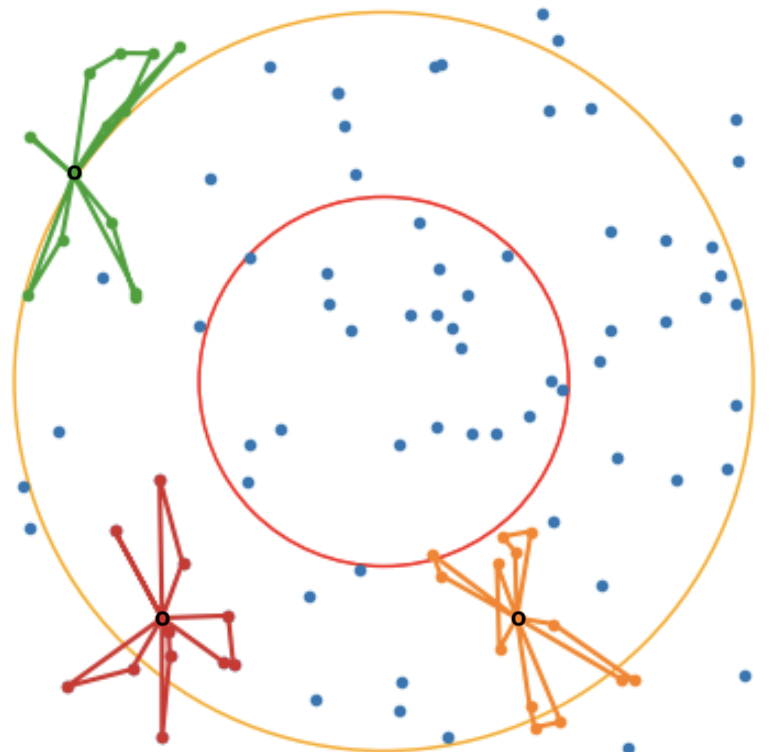}
   \caption{Warehouse distribution in the case where population distribution is uniform, and traffic zones are present}
   \label{fig:case3}
\end{figure}

\begin{enumerate}
    \item Uniform traffic and population distribution (Fig. \ref{fig:case4}): This is the base case, to observe how other cases orient themselves compared to it. In an area with uniform population distribution and travel times, the warehouses do not have any pre-determined parameter to base a preference upon. 
\item Variable traffic and gaussian population distribution (Fig. \ref{fig:case1}): Here there are two contrasting parameters which counter each other. The gaussian population distribution pushes the warehouse to be within the central region, while the higher travel times counter this effect pushing them further out. This case is closest to what happens in reality in most cities, where the central region of the city has a high population distributions, but higher traffic than the outer regions. 
\item Uniform traffic and gaussian population distribution (Fig. \ref{fig:case2}): In this case, the warehouses tend to be within the central region, considering it has the highest population distribution with the same amount of travel times as the outer regions, giving it a considerable bias. This case manifests itself in cities that have excellent traffic management in central or downtown regions, prompting the warehouses to be located within the busiest areas. 
\item Variable traffic and uniform population distribution (Fig. \ref{fig:case3}): In this case, the warehouses tend to be pushed towards the outer areas of the region, where traffic times are lower even though the population density is the same.
\end{enumerate}


After running the experiment on python, using randomly generated points for customer locations as the dataset with sizes of 50, 100 and 200 in the variable traffic zone case and 100, 200 and 400 in the uniform travel time case, we have the following observations. For the both the uniform travel time and variable traffic zone cases, we observe an increase in time on average by a factor of 7 on doubling the sample size. 

From the various scenarios observed in this experiment, we get a strong idea of how a quick delivery app would base its decisions to install dark stores within a given city they wish to service, with the goal of profitability in mind. These apps generally start with servicing one area, before installing further warehouses in other parts of the city, and this is exactly the same approach which we have followed as well.

\textbf{Declaration of generative AI and AI-assisted technologies in the writing process}

During the preparation of this work the author(s) used ChatGPT in order to improve readability and cogency in text. After using this tool/service, the author(s) reviewed and edited the content as needed and take(s) full responsibility for the content of the publication.

\bibliographystyle{plain}
\bibliography{main}

\end{document}